\magnification\magstep1
\widowpenalty10000
\font\smc=cmcsc10
\font\sans=cmss10
  %%% want 7
\font\bigbf=cmbx12
\def\ce{\centerline}

\input amssym.def
\newsymbol\Vdash 130D
\newsymbol\nVdash 2331
\newsymbol\backepsilon 237F
\def\a{\alpha}
\def\vphi{\varphi}
\def\phi{\vphi}  %%%%%%%%%%%%
\def\B{\hbox{\sans{I$\!$B}}}
\def\P{\hbox{\sans{I$\!$P}}}

\def\cD{{\cal D}}
\def\cN{{\cal N}}
\def\cT{{\cal T}}

\def\nmodels{\not\models}     %%%%
 %%%%
\def\fo{\Vdash}
\def\nfo{\nVdash}

\def\implies{\to}
\def\contains{\!\backepsilon}
\def\and{\wedge}
\def\iff{\leftrightarrow}

\def\ese{\emptyset}
\def\eps{\,\epsilon\,}
\def\veps{\,\varepsilon\,}
\def\vepsa{\,\varepsilon_\a\,}

\def\si#1{\buildrel #1 \over \sim}

\def\pf{\medskip Proof:\ }
\def\Lemma{\medskip{\bf Lemma: }}
\def\eop{\hfill\medskip}

\ \bigskip %%%%

\ce{\bigbf Foundations for abstract forcing}
\medskip
\ce{Peter M.\ Johnson \qquad {\tt peterj@ufba.br}} %(V6, Sept.\ 2007)} 
\medskip
\ce{Universidade Federal da Bahia, Salvador, Brazil} 

\bigskip
{\bf Abstract:}
The foundations of forcing theory are reworked to streamline the presentation
and to show how the most basic results are applicable in very general contexts. 
\bigskip
\ce{\smc Introduction}
\medskip

After Cohen [Co] invented the method of forcing to settle long-open independence results 
in set theory, techniques evolved to become simpler as well as more powerful.  
Our focus is on simplifying foundational aspects of forcing, in its most
elementary form, that seem to have been relatively neglected.   
As will be seen, forcing becomes particularly easy to understand when separated 
into distinct but interacting techniques, involving posets, generic filters, 
Boolean algebras, and (for set theory) a weak form of the membership relation.
A new general-purpose forcing method, with a presumably much wider range of
applications, arises when aspects related to the construction of models 
for set theory are stripped away and the role of names reexamined.
Well-known complications involving the set membership relation $\in$ do not arise 
if methods requiring well-foundedness are abandoned and simpler ones found.
Although forcing models are studied via objects called names, constructed in a 
standard way (with minor variations) from the poset, details about names will 
turn out to be irrelevant in the earlier foundational stages studied below.
Could new ways of encoding names, with new relations and interpretations of them, 
yield innovative applications of forcing in areas distant from set theory?

An overview of the most relevant aspects seems appropriate.  Forcing 
has been regarded mainly as a method to construct models of set theory tailored to have 
a great variety of properties, resolving purely set-theoretic questions such as those 
concerning possible values of powers of cardinals, or more `applied' ones in areas 
that include infinite combinatorics, topology, measure theory, and proof theory.  
Very briefly, one approach starts with a fixed countable transitive 
model $(M, \in)$ of ZFC which gives rise, nonconstructively, to other more 
interesting models.  It is convenient to work with such an $M$ but
its existence is an unecessarily strong assumption,
and several alternate approaches are possible (e.g., see [Ku]).
The most crucial aspect in applications is to define a poset, often with
great ingenuity, so that it `forces' the new models obtained using it to have 
the properties aimed at.
These properties are discussed in a forcing language, using constants intended 
to represent (nonfaithfully) names of the intended new sets.  Statements are interpreted 
in a way that depends on the choice, outside $M$, of a `generic' filter in the poset
and on the form in which names are represented concretely as certain sets in $M$.
This approach already appears in mature form in Shoenfield [Sh], whose 
presentation incorporated important simplifications.  For historical details, see [Mo].
Since the basic ideas of forcing lead to spectacular applications, some easily 
accessible, they deserve to be made as simple and transparent as possible.

\goodbreak %%%%
\bigskip
\ce{\smc Order-theoretic preliminaries}
\medskip

Order-theoretic concepts usually refer to the fixed poset $(\P, \leq)$,
and variables such as $p,q,r$ range over the elements of $\P$. 
Notation $X^{\uparrow} = \{p : (\exists x \in X)(x \leq p)\}$ and $X^\downarrow$ 
is used for up- and down-closures of subsets $X$ of $\P$.
Optionally, one can use topological concepts, relative to the topology on $\P$
formed by the down-closed sets.  For example, a subset $D$ is {\sl dense} if 
$D^{\uparrow} = \P$.

Instead of working with arbitrary posets, one can use Boolean algebras.
This permits techniques of a more algebraic nature and more concise language.
In practice, to construct a model tailored to have certain properties, the 
preferred method is search for an appropriate poset $\P$, which is then embedded into a
complete Boolean algebra.  
One way is to start with the unary operation which takes each subset $X$ of $\P$  
to a down-closed set $X' = \P \backslash X^\uparrow$.
This yields a closure operator $X \to X''$ on $\P$, whose closed elements
($X = X''$) form a family, here called $\B(\P)$.   
As usual, the family is closed under arbitrary intersections, making it
a complete lattice.
The lattice is complemented and is in fact a Boolean algebra, by a result of Byrne 
[By], as it satisfies Byrne's axioms: $\and$ is a semilattice operation, $0 \neq 0'$, 
and, most notably, $X \and Y' = 0$ iff $X \and Y = X$.
Alternatively, $\B(\P)$ consists of the regular open sets of the above topology on $\P$.

In this paragraph, we make the mild assumption (separativity) that the sets 
$p^\downarrow$ $(p \in \P)$ lie in $\B(\P)$.  Then $\P$ has 
a certain kind of universal embedding in a complete Boolean algebra. 
One good source for relevant details is [Kub].
We remark that $\P \to \B(\P)$, $\,p \mapsto p^\downarrow$ is such an 
embedding.  As this will not be needed, the details are omitted.

In conventional applications of forcing, the poset $(\P,\leq)$ lies within some model 
$(M,\in)$ for set theory, and those elements of $\B(\P)$ that lie in $M$  
clearly form a subalgebra $\B_M(\P)$, with sup and inf
restricted to subfamilies definable in $M$.  This relativization of $\B(\P)$, 
viewed within the model, is a complete Boolean algebra. 
Moreover, the forcing map (see below) which takes a statement $\phi$ to
$\{ p \in \P : p \fo \phi \}$ will be seen to map into $\B_M(\P)$.

A filter $G$ is a non-empty up-closed set for which each pair of elements of $G$
has  a lower bound in $G$.
The dense subsets of $\P$ are clearly those of the form $X \cup X'$, and
no filter can contain points of both $X$ and $X'$.  
Relative to any collection $\cD$ of dense subsets, there is a notion of $\cD$-generic 
filters $G$, by which it is meant that $G$ has non-empty intersection with each 
$D \in \cD$.  A well-known idea (see [Sh] or [Ku, VII, 2.3]) produces an abundance of
$\cD$-generic filters, enough so that each point $p$ lies in such filters,
whenever $\cD$ is countable.  

\bigskip
\ce{\smc Abstract forcing}
\medskip

A {\sl forcing language} will be taken here to be just a first-order language, say
with relation symbols but not terms, used to make statements involving
constants called `names', which are used in place of free variables.
The choice of basic symbols for logical conectives and quantifiers 
(say $\neg,\and,\exists$) matters little in our classical approach.
If desired, one could rephrase results in terms of a version of the 
Lindenbaum algebra, regarded as the Boolean algebra (with some extra structure)
naturally induced on the statements of the forcing language modulo logical equivalence.

Variables such as $\tau,\sigma,\pi$ always range over names, and may be quantified.
The abstract approach to names requires nothing more than a function from atomic 
statements $\phi$ to values $[[\phi]]$ in the Boolean algebra $\B(\P)$. 
In set-theoretical applications, the class of names and the relation $p \in 
[[\sigma \in \tau]]$ should be definable by predicates in the base model $(M,\in)$.
Extensionality of the new models is obtained by imposing a condition $(E)$, treated below.

Any evaluation of atomic statements extends, in the obvious way, to one for 
all statements in the language.
In particular, $[[(\exists \tau)\phi(\tau)]]$ is the sup in $\B(\P)$ 
of the elements $[[\phi(\tau)]]$ as $\tau$ ranges over all names.
Without assuming more about the structure of names,
this sup is not usually the union, which could be thought of as a defect when 
working with individual points in $\P$.  Generic filters provide the remedy.

The countable collection $\cD$ used to define `generic' matters little, so long as it 
contains all $X \cup X'$ for which $X$ is either of the form $[[\phi]]$ or (letting
one of the names used in $\phi$ vary freely) of the form $\cup_\tau [[\phi(\tau)]]$.  
This restriction on $\cD$ is essential for the next results, but some class 
$\cD$ more easily definable than the minimal one is customarily used.  
To ensure that $\cD$ is countable, the approach for set-theoretic applications makes the 
strong assumption that $M$ (viewed from outside) is countable, so that only a countable 
number of dense subsets of $\P$ are definable in $M$.  In more general contexts, it would 
suffice to assume at most enumerably many names and relation symbols (hence also statements).

For each generic filter $G$, consider the set $\cT_G$ of statements $\phi$ for which 
$[[\phi]]$ intersects $G$ nontrivially; i.e., $(\exists p \in G)(p \in [[\phi]])$.
These $\phi$ are, in some sense `true relative to $G$', and we write $G \models \phi$,
while statements not in $\cT_G$ are `false for $G$'. 
This assignment of truth values respects operations such as $\and,\neg$.
To obtain models for $\cT_G$, quantifiers must also be respected. 
Handling $\exists$ yields at once the perhaps more surprising result for $\forall$. 

\Lemma
For statements $\phi(\tau)$, with $\tau$ varying over names, and $G$ a generic filter, 
\item{(a)} 
$\ G \models \exists \tau\phi(\tau)$ \ iff \ for some $\tau$, $G \models \phi(\tau)$.
\item{(b)} 
$\ G \models \forall \tau\phi(\tau)$ \ iff \ for all $\tau$, $G \models \phi(\tau)$.

\pf
For this, it suffices to show that, given 
any family of statements of the form $\phi(\tau)$ with $\tau$ ranging over names, and 
letting $X$ be the union of the corresponding subsets $[[\phi(\tau)]]$, 
then $[[(\exists \tau)\phi(\tau)]]$ is disjoint from $G$ if $X$ is.  
This holds because here $G$ must contain some point of $X'$, hence none in $X''$.
This last set is just $[[(\exists \tau)\phi(\tau)]]$.
\eop

It should not fail to be mentioned that, for names of the form used in set theory, 
part (a) of the lemma can be strengthened to say that some suitable name $\tau$ 
depends only on $\phi$, not on $G$.  See for example [Kub, Th.~9.2].
(The treatment in [Ku, VII, 8.2] is incomplete.) 
Although statements $\phi$ are evaluated in the base model by using the 
Boolean algebra, the usual aim is to use a generic filter $G$ to construct 
a model $M[G]$ with $M[G] \models \phi$ iff $G \models \phi$.
Note that the previous lemma implies at once that the forcing language has such 
a model: for now, $M[G]$ will be taken to be the set $\cN$ of names, endowed via the 
atomic statements in $\cT_G$ with the relations corresponding to the relation symbols.
Later, one may prefer to factor out $\cN$ by some equivalence relation, 
or use an appropriate image of $\cN$, in order to interpret more naturally
symbols such as $=$ (which should be equality) and $\in$ (membership).
  
Each statement $\phi$ in the forcing language yields a forcing relation, say
as a unary relation on $\P$, as follows.
For points $p$ of $\P$, we say that $p$ {\sl forces} $\phi$, and write $p \fo \phi$, 
when $\phi$ is true in every model $M[G]$ with $G$ generic (relative to $\cD$) and
$G \contains p$.

There are two especially fundamental results concering forcing.  First, it is
definable in a way that does not mention generic filters: for each statement $\phi$, 
some predicate in the base model determines which $p \in \P$ force $\phi$.
The second, the Truth Lemma, is that a statement $\phi$ is valid in a 
model $M[G]$ exactly when some $p$ in $G$ forces $\phi$.
Another perhaps noteworthy observation is that the forcing relation does not depend 
on the choice (restricted as above) of the dense sets $\cD$ used to define generic filters. 
All these results follow immediately from:
\Lemma  $p \fo \phi$ iff $p \in [[\phi]]$.
\pf If $p \in [[\phi]]$ then, trivially, $p \fo \phi$.  Conversely, when $p \notin [[\phi]]$ 
some $q \leq p$ lies in $[[\phi]]'$, and some generic filter $G$ contains $q$.   
Now $M[G] \nmodels \phi$ and $p \in G$, giving $p \nfo \phi$.

\bigskip
\ce{\smc Interlude on non-extensional and non-well-founded theories}
\medskip

In theories and models, the symbol $=$ is almost universally taken to mean equality.
Following a principle held by Leibniz, objects which cannot be distinguished by 
any relevant property should be regarded as equal.  
The only reason not to identify would be to leave open the possibility of later more 
refined properties. In models of set theory, or within the language itself, $=$ can be 
treated as being defined from the basic relation $\in$, via an even 
stronger principle, extensionality.  This idea, that sets with the same 
members are equal, seems so central to the concept of what sets are (platonically) 
that it has rarely been examined critically.  When $\in$ is the only basic relation 
present, Leibniz extensionality, or quasi-extensionality (with an equivalence 
relation in place of $=$) is:
$$(E) \qquad (\forall \pi)(\pi \in \sigma_1 \iff \pi \in \sigma_2) \implies 
(\forall \tau)(\sigma_1 \in \tau \iff \sigma_2 \in \tau).$$
However, set-theoretic forcing tends to produce even laxer models that need 
further adjustments in order to satisfy (E). 
Another principle, that sets should be well-founded, is usually imposed here,
but this a rare example where such ideas only serve to complicate matters.  
Now that there is sufficient motivation to examine non-extensionality, 
in conjunction with non-well-foundedness, we take a broad view of related issues.

The area of mathematical foundations is rich with history and 
interesting issues, enough to merit continuing scholarly attention.
Attitudes to certain topics, for example foundations for set theory, may change 
for reasons worth elucidating, detailing how groups with different agenda, using 
different language, influence each other or fail to do so.   For example, it is
extremely convenient to work with well-founded sets, using transfinite
induction on rank.  This may explain in part why non-well-founded set theories 
remained a little-studied curiosity in the shade of ZFC until there was a clear 
need for them, as tools useful for studying several problems in computer science.
For similar reasons, non-extensionality (and intensionality) have in recent years 
become common topics in computer science, whose influence may in time spread wider.

Despite earlier work in the field, notably by Boffa (in many articles), the use 
of non-well-founded sets became commonplace thanks to the timely and influential 
work by Aczel [A]. 
Not surprisingly, the sort of non-well-founded theories considered are overwhelmingly
those based on Aczel's anti-foundation axiom AFA, which is too restrictive for our taste.  
Our view is that if objects and a binary relation called $\in$ are somehow given,
and objects with the same `members' are equal, then this {\sl is} a system of sets, 
as long as we do not need to examine objects internally (using more than $\in$) to see 
what they are `really' made of.  Does anyone ask what the empty set is made of?  As in
other areas of mathematics, only abstract structural properties matter.  Following 
Boffa rather than Aczel (and others), we regard it as desirable to permit sets with
different elements to have isomorphic {\sl internal} $\in$-structures, thereby for
example allowing more than one set of the form $x = \{x\}$.  

Many equivalent ways have been used to axiomatize the system known as ZF.
The exact form of the axioms becomes important when studying systems which omit or 
modify axioms, and our attention is on Extensionality together with Foundation.  
One standard source treating axioms for ZF and other set theories is [F],
but it makes no claim to completeness and mentions surprisingly little 
about the particular topics we wish to focus on. 
For non-extensional theories, the following studies, whose reviews we consulted, seem 
to be among the most mathematically relevant.
Hinnion [Hi] introduces methods which, among other things, 
simplify earlier work by Gandy [Ga], Scott [Sc], and others. 
Briefly, the effects of omitting Extensionality are as follows.
Scott showed that the resulting system ZF${}^{\neq}$ is distinctly weaker than ZF,
while Gandy and H.\ Friedman [Fr] showed that their systems suffer no significant change.
The reason is that Scott and Gandy used the usual form of Replacement, involving a 
functional relation of the form $(\forall x)(\exists! y)\vphi(x,y)$, but Gandy
compensated for this weakness by introducing a set formation operator 
$\lambda x \,A(x)$ (i.e., $\{x \mid A(x)\}$). 
The solution we prefer is Friedman's, who used the Collection axiom (not even
mentioned in [F]), which strengthens Replacement when Foundation is not assumed.
Collection is: $[(\forall x)(\exists y)\vphi(x,y)] \to [(\forall X)(\exists Y) 
(\forall x \in X)(\exists y \in Y)\vphi(x,y)]$. 
A good exposition of Friedman's later improvements on [Fr] appears in [Kr, \S5].
The conclusion is that, even without Foundation, extensional systems can easily
be recoved from non-extensional ones of the right form.  

Friedman's solution is internal to the set theory ($\in$ can be defined from a 
weaker non-extensional $\veps$, assuming some set-theoretic axioms, but not 
Foundation), whereas our solution assumes no axioms, but is carried out within 
the Boolean algebra rather than within the forcing language.
It can also be carried out internally in the presence of axioms sufficient 
to support a definition by transfinite induction.
While Friedman recovers extensionality by identifying objects as much as possible, 
using an approach based on ideas of Aczel, where the preferred bisimulations are 
maximal, we identify objects only when absolutely necessary to obtain extensionality.  
The method is so simple that it may have been rediscovered several times, but it 
does not seems to be widely known, as intensive searches failed to locate a 
refererence in literature or reviews accessible to us. 
Aczel [A] takes a different approach when studying the least bisimulation, 
a topic treated earlier by Hinnion [Hi], presumably in a roughly similar way.
However, the method we prefer can be found in two recent preprints -- 
[BMW, Prop.~6] and [Fi, \S6].
In the second, the idea is used in the process of making an interesting 
connection between forcing and modal logic.
 
The solution now given produces the smallest quasi-extensional relation $\in$ (with 
an equivalence relation in place of equality) generated from an arbitrary binary 
relation $\veps$ on a set.  The relation $\in$ is well-founded precisely when $\veps$ is.
The language will temporarily have many binary relations: 
$\in, \veps$, $\sim$ and a family $\si\a$, where $\a$ always ranges over ordinals.
For abstract forcing, the universe is the class of names, and only ordinals $\a < |\P|$ 
are needed.  For more general use, name-free notation will be used.

To start, $\si0$ is equality, and $y_1 \si{\a+1} y_2$ means 
$(\forall x_1)( x_1 \veps y_1 \implies 
(\exists x_2) (x_1 \si\a x_2 \and x_2 \veps y_2)) \and
(\forall x_2)( x_2 \veps y_2 \implies 
(\exists x_1) (x_1 \si\a x_2 \and x_1 \veps y_1)).$
For limit ordinals $\lambda$, $\si\lambda$ is the limit (union) of the earlier relations.
It is not hard to see that the $(\si\a)_\a$ form an increasing family of 
equivalence relations between names, and $\sim$ is defined to be its limit.
Finally, $x \in y$ is defined to mean $\exists x' (x \sim x' \and x' \veps y)$.
It may be clearer to imagine the transition from $\veps$ to $\in$ in stages,
defining $x \vepsa y$ iff there exist $x' \si\a x$
and $y' \si\a y$ with $x' \veps y'$.
The earlier extensionality condition $(E)$ clearly holds.
The equivalence relation $\sim$ (which is factored out to form models with $=$) and 
the relation $\in$, are definable in terms of each other, in the presence of $\veps$.  
Thus $x \veps y \implies x \in y$ and $\veps$ is a simulation for $\sim$, by which we 
mean that $x \veps y \and y' \sim y \implies (\exists x' \sim x)( x' \veps y')$.
This gives a clear idea how to construct all possible $\veps$ from a 
quasi-extensional relation $\in$ (with an equivalence relation $\equiv$
in place of $=$), not guaranteeing here that
$\equiv$ will be the smallest relation $\sim$ compatible with $\veps$.

\bigskip
\ce{\smc 
Forcing with names and Boolean-valued models}
\medskip

Returning to forcing, with the enriched language now using names, 
formulas will be assigned values in a complete Boolean algebra $\B$,
starting with atomic formulas and extending in the obvious way.
Everything is done within a model $(M, \in)$ for set theory, in ways 
definable by predicates in the forcing language, using a given 
valuation of $\veps$ on the class of names, with values in the Boolean algebra
$\B = \B_M(\P)$ constructed from a fixed poset $\P$ in the model.  
If $M$ is countable, we can also concentrate on a model of the forcing language 
obtained from a generic filter $G$.  This initially gives a non-extensional binary 
relation $\varepsilon_G$ on the names, which generates a quasi-extensional relation 
$\in_G$ in the way described in the previous section.   As well-foundedness is not 
assumed, it no longer seems compelling to try to obtain quotient models $(M[G], \in)$ 
that are standard.

The construction using Boolean valuations starts by defining
$[[\tau_1 \si0 \tau_2]]$ to be 1 (the set $\P$) if the names $\tau_1$ 
and $\tau_2$ are identical, and 0 ($\ese$) otherwise. Whatever $G$ is, this just 
gives the equality relation on names.  Other definitions are as expected.
Thus, for limit ordinals $\lambda$, 
$[[\tau \si\lambda \tau']] = \bigvee_{\a < \lambda} [[\tau \si\a \tau']]$.
Care is required here, as it is conceivable that a generic filter $G$ could be 
disjoint from the sets $[[\tau \si\a \tau']]$ $(\a < \lambda)$ but not from 
$[[\tau \si\lambda \tau']]$.
In set-theoretic models, assuming enough axioms to allow quantification 
over ordinals, an earlier lemma shows that this problem does not arise.
In general, a suitable relation $\sim$, no longer claimed to be the smallest one,
will be produced.  This follows from $[[\tau \si\lambda \tau']] = 
\bigvee_{\a < \lambda}[[\tau \si\a \tau']] \leq [[\sigma \veps \tau \implies 
(\exists \sigma' \veps \tau')(\sigma \si\lambda \sigma')]]$. 
 
Set-theoretic axioms are usually stated in terms of $\in$,
but stronger forms using $\veps$ may be more convenient.
The axioms of ZFC are especially well-suited for constructing models 
that satisfy enough of the axioms (a finite number) to be useful.
Names and Boolean valuations can be built in transfinite stages from initial 
data, in ways ensuring that these axioms then hold in the models $(M[G], \in_G)$.
Now that problems concerning non-extensionality and non-well-foundedness have 
been adequately resolved, remaining details can be inferred from many sources.
Only a little will be sketched here.

To construct names of one of the conventional forms,
one can start with $\cN_0 = \ese$ (say), define $\cN_{\a+1}$ to be the union of 
$\cN_\a$ with the set of functions $\cN_\a \to \B$ that lie in $M$, and let
$\cN_\lambda$ ($\lambda$ a limit ordinal) be the obvious union.
Thus each name is `created' at a successor ordinal in $M$.
For names $\sigma, \tau$, the value $[[\sigma \veps \tau]]$ is defined to be
$\tau(\sigma)$ if $\sigma$ is created before $\tau$, and is 0 otherwise.

As an illustration, the power set axiom will be examined.
Write $\sigma' \leq \sigma$ if, for all names $\pi$, $[[\pi \veps \sigma']] \leq 
[[\pi \veps \sigma]]$.  This provides a sufficiently large supply of names for subsets.
Given $\sigma \in \cN_\a$ one name $\tau: \cN_\a \to \B$ with
$[[\forall \sigma'(\sigma' \in \tau \iff \sigma' \subset \sigma)]] = 1$     
is as follows:  for all $\sigma' \in \cN_\a$, $\tau(\sigma') = 1 \in \B$ when 
$\sigma' \leq \sigma$, while $\tau(\sigma') = 0$ otherwise.

{\bigskip
\parskip3pt
\ce{\smc References}
\frenchspacing
\medskip

[A] Aczel, P.: {\sl Non-Well-Founded Sets},  CSLI Lecture Notes {\bf 14}, Stanford, 1988.

[BMW] Bab, S.; Marr, B.; Wieczorek, T.: {\sl $\eps$-Style (of) Semantics.
An alternative to set-theoretic modelling}, preprint, 31pp.

[Bo] Boffa, M.: {\sl Forcing et n\'egation de l'axiome de fondement},
Acad. Roy. Belgique, Mem. Cl. Sci., Coll. 8$^\circ$, II. S\'er., Tome {\bf XL}, No.7 (1972), 53pp.
% [Bo2] others by Boffa?

[By] Byrne, L.: {\sl Two brief formulations of Boolean algebra}, 
Bull. Amer. Math. Soc. {\bf 52} (1946), 269--272.

[Co] Cohen, P. J.: {\sl The independence of the continuum hypothesis, I},
Proc. Nat. Acad. Sci. USA {\bf 50} (1963), 1143--1148.
% II  51 (1964) 105--110.

[F] Fraenkel, A., Bar-Hillel, Y., Levy, A.: {\sl Foundations of Set Theory}, 
$2^{\rm nd}$ ed., vii+404 pp., North-Holland, 1973.

[Fi] Fitting, M.: {\sl Intensional Logic -- Beyond First Order}, preprint, 22 pp., c.\ 2003.

[Fr]  Friedman, H.: {\sl The consistency of classical set theory relative to 
a set theory with intuitionistic logic},
J. Symbolic Logic {\bf 38} (1979), 315--319.

[Ga]  Gandy, R.: {\sl On the axiom of extensionality, II},
J. Symbolic Logic {\bf 24} (1959), 287--300.
% I: {\bf 21} (1956), 36--48; II,  % idem, 

[Kr]  Krivine, J.-L.: {\sl Typed lambda-calculus in classical Zermelo-Fraenkel
 set theory}, Arch. Math. Log., {\bf 40}(3) (2001), 189--205.  % no.3

[Kub] Kubis, W.: {\sl Forcing}, preprint (1999), 21 pp.

[Ku] Kunen, K.: {\sl Set Theory}, xvi+313 pp., North-Holland, 1980.
 
% [L] Luckhardt, H.: {\sl Extensional G\"odel Functional Interpretation}, 
% % Lecture Notes in Mathematics
% LNM {\sl 306}, Springer, Berlin, 1973. 
%
[Mo] Moore, G. H.: {\sl The orgins of forcing}, 
% Logic Colloquium ´86 (Hull),
Stud. Logic Found. Math. {\bf 124} (1988), 143--173.
%
% [Po] Pollard, S.: {\sl A Strengthening of Scott's ZF${}^{\neq}$ Result},
% Notre Dame J. of Formal Logic {\bf 31}(3) (1990), 369--370.  % no.3
%
% [Ro] Robinson, A.: {\sl On the independence of the Axioms of Definiteness 
% (Axiome Der Bestimmtheit)}, J. Symbolic Logic, {\bf 42}(1939), 69--72. % June

[Sc] Scott, D.: {\sl More on the Axiom of Extensionality}, pp. 115-131 in {\sl
Essays on the Foundations of Mathematics}, ed. Bar-Hillel, Y. et al., 
Hebrew University, Jerusalem, 1966. 

[Sh] Shoenfield, J. R.: {\sl Unramified Forcing}, in {\sl Proc. Symp. Pure Math.}
 {\bf 13}(1971), 357--381.

% [Sa] Sato Kentaro, {\sl Forcing under Anti-Foundation Axiom: An expression of the 
% stalks}, Math. Log. Quart. {\bf 52} (2006), 295--314.

}

\bye